\documentclass[12pt,secthm]{elsart}
\usepackage{amsmath,amssymb,amscd}
\input xy
\xyoption{all}
\makeatletter
\renewenvironment{pf}%
  {\par\addvspace{\@bls \@plus 0.5\@bls \@minus 0.1\@bls}\noindent
     {\bfseries\Elproofname}\enspace\ignorespaces}%
       {\par\addvspace{\@bls \@plus 0.5\@bls \@minus 0.1\@bls}}
       \def\Elproofname{%
        Proof:
        }
\makeatother

\theoremstyle{plain}
\theoremstyle{definition}

\DeclareMathOperator{\Tr}{{Tr}}

\DeclareMathOperator{\End}{{End}}

\def\lra{\longrightarrow}

\def\Z{{\mathbb Z}}
\def\Q{{\mathbb Q}}

\def\P{{\mathbb P}}

\def\M{{\mathcal M}}

\begin{document}
\begin{frontmatter}
\title{A family of Prym-Tyurin varieties of \\exponent 3}

\author[lange]{H. Lange} \ead{lange@mi.uni-erlangen.de} ,
\author[recillas]{S. Recillas} \ead{sevin@matmor.unam.mx} ,
\author[rojas]{A.M. Rojas} \ead{anita@matmor.unam.mx}

\address[lange]{Mathematisches Institut,
              Universit\"at Erlangen-N\"urnberg, Germany}
\address[recillas]{Instituto de Matematicas,
              UNAM Campus Morelia and CIMAT, Guanajuato, M\'exico}
              
\address[rojas]{Instituto de Matematicas,
              UNAM Campus Morelia, M\'exico}

\thanks{Supported by DAAD, Conacyt 40033-F and FONDECYT No.3040066}

\begin{abstract}
We investigate a family of correspondences associated to \'etale coverings of degree 3 of hyperelliptic curves.
They lead to Prym-Tyurin varieties of exponent 3. We identify these varieties and derive some consequences.
\end{abstract}

\end{frontmatter}


\maketitle

\section{Introduction}

A {\it correspondence} on a smooth projective curve $C$ is by definition a divisor $D$ on the product $C \times C$. Any 
correspondence $D$ on $C$ induces an endomorphism $\gamma_D$ of the Jacobian $JC$. Conversely, for every 
endomorphism $\gamma \in \End C$ there is a correspondence $D$ on $C$ such that $\gamma = \gamma_D$. However 
for most correspondences $D$ which occur in the literature, $\gamma_D$ is a multiple $d\cdot 1_{JC}$ of the 
identity of the Jacobian. In particular this is the case for any correspondence of a general curve. 
It were mainly these correspondences, called {\it of valency d}, which were studied 
by the classical Italian geometers (see e.g. \cite{sev}).\\
At the beginning of the 1970's A. Tyurin suggested the investigation of another class of correspondences, 
namely effective symmetric correspondences $D$ without fixed point on $C$ such that $\gamma_D$ satisfies an equation
$$
\gamma_D^2 + (e-2)\gamma_D -(e-1) = 0.
$$ 
For these correspondences $P = im(\gamma_D - 1_{JC})$ is a {\it Prym-Tyurin variety} of exponent $e$, meaning that the restriction of the 
canonical polarization of $JC$ to $P$ is the $e$-fold of a principal polarization on $P$. Jacobians are Prym-Tyurin 
varieties of exponent 1, Prym varieties associated to \'etale double coverings are Prym-Tyurin varieties of exponent
2. On the other hand, it is not difficult to show (see \cite{lb}, Corollary 12.2.4) that every principally polarized 
abelian variety is a Prym-Tyurin of some high exponent. However it seems not so easy to construct Prym-Tyurin varieties 
of low exponent $\geq 3$. First examples (associated to Fano threefolds) were investigated by Tyurin (see \cite{tyu}).
Other examples (associated to Weil groups of certain Lie algebras) were constructed by Kanev (see \cite{kan}).

It is the aim of this paper to study the following correspondence: Let $C$ be a hyperelliptic curve of genus $g \geq 3$  
and $f: \tilde{C} \lra C$ an \'etale threefold covering. Consider the following curve in the symmetric product 
$\tilde{C}^{(2)}$
$$
X = \{p \in \tilde{C}^{(2)}\;|\; f^{(2)}(p) \in g^1_2 \}
$$
Let $\iota$ denote the hyperelliptic involution of $C$ and for $x \in C$ write $f^{-1}(x) = \{x_1, x_2,x_3\}$, 
$f^{-1}(\iota x) = \{y_1,y_2,y_3\}$ and moreover for abbreviation $P_{ij} = x_i + y_j \in \tilde{C}^{(2)}$. Then 
the symmetric $(2,2)$-correspondence $D$ on $X$ is defined by
$$
D = \{(P_{ij},P_{kl}) \in X \times X \;|\; i=k \, \mbox{and} \; j \not= l \; \mbox{or} \; i \not= k \;\mbox{and} \; j=l \}.
$$ 
We show in section 2 under the hypothesis that $X$ is smooth and irreducible, that $P = im (\gamma_D - 1_{JX})$
is a Prym-Tyurin variety of exponent 3. In section 3 we realise a $(2g-1)$-dimensional family of pairs $(C,f)$ such that
$X$ is smooth and irreducible. In fact, the Galois group $G$ of the Galois extension 
$Y \lra \P^1$ of $\tilde{C} \lra \P^1$ is necessarily isomorphic to 
$S_3 \times S_3 \subset S_6$. In section 5 we compute the dimensions of the Jacobians and Prym varieties relevant
to this situation. The main result of this section is Theorem 5.3 which 
says that there are two trigonal curves $X_1$ and $X_2$ associated to subgroups of $G$ such that $P$ is 
canonically isomorphic as a principally polarized abelian variety to the product of Jacobians $JX_1 \times JX_2$. The trigonal covers of $X_1$ and $X_2$ have 
disjoint ramification locus and $X$ is their common fibre product over $\P^1$. As a consequence of this and the moduli
considerations of section 4 we obtain the following consequence which seems of interest to us and for which we could 
not find a different proof:\\

\noindent
{\bf Corollary} (of Theorems 4.1 and 5.2): {\it 
Let $X_1$ and $X_2$ be trigonal curves with simple ramification and disjoint branching and let $X$ denote their fibre product
over $\P^1$ with projections $f_i: X \lra X_i$. Then $X$ is a smooth projective curve and 
the map $f_1^* + f_2^*: JX_1 \times JX_2 \lra JX$ is an embedding.}\\

Finally in section 6 we study the Abel-Prym map of the Prym-Tyurin variety $P$ defined as the composition of the 
Abel map $X \lra JX$ and the projection $JX \lra P$. The main result is the following\\

\noindent
{\bf Theorem}
{\it For $g \geq 6$ the Abel-Prym map $\beta_P: X \lra P$ is an embedding.}\\

As a consequence we obtain that the cohomology class of three times the canonical product polarization of $JX_1 \times JX_2$
is represented by a smooth irreducible curve (see Corollary 6.5). 

\section{Construction of the Prym-Tyurin varieties}\label{S:constructionpt}

Let $C$ be an hyperelliptic curve of genus $g \geq 3$, $i:C \to C$ the hyperelliptic involution and $h:C \to \P^1$
the map given by the $g_2^1$.
Let $f:\tilde{C} \to C$ be an \'etale morphism of degree $n$ from a projective smooth irreducible curve $\tilde{C}$.

Let us define a new curve $X$ by the following cartesian diagram

$$
\xymatrix{
  X:=(f^{(2)})^{-1}(g_2^1) \ar[d]_{\pi=f^{(2)}|_X} \ar@{^{(}->}[r]
                & \tilde{C}^{(2)} \ar[d]^{f^{(2)}}  \\
  \P^1\cong g_2^1  \ar@{^{(}->}[r]   & C^{(2)}             }  \eqno(2.1)
$$    

\noindent where $f^{(2)}: C^{(2)} \lra \tilde{C}^{(2)}$ denotes the second symmetric product of $f$. 
Observe that $\pi$ is of degree $n^2$.

For the rest of this section let us assume that the curve $X$ is {\it smooth and irreducible}. (In the next section we will see
that there exist \'etale coverings $f$ such that this is the case.) Under this hypothesis we
define a correspondence $D$ on $X$. For this consider the canonical $2:1$-map $\lambda : \tilde{C}^2 \lra \tilde{C}^{(2)}$ from the 
cartesian to the symmetric product of $\tilde{C}$ and denote 
$$
\tilde{X} := \lambda^{-1}(X).
$$ 
Let $p_1: \tilde{X} \lra \tilde{C}$
denote the projection onto the first factor, where $\tilde{X}$ is considered as a curve in
$\tilde{C}^2$. Then
$$
\tilde{D} := \{(a,b) \in \tilde{X} \times \tilde{X} \,\, | \,\, p_1(a) = p_1(b) \}
$$
with reduced subscheme structure is an effective divisor on $\tilde{X}^2$ containing the diagonal $\tilde{\Delta}$.  
Denote
$$
Y := \tilde{D} - \tilde{\Delta}
$$
The divisor
$$
D := (\lambda \times \lambda)_*(Y)
$$
is an effective symmetric correspondence of $X$ of bidegree $(2n-2,2n-2)$.  

In order to describe this correspondence set-theoretically we fix some notation. Given $z \in \P^1$, let 
$$
h^{-1}(z)=x+ix
$$ and
$$
f^{-1}(x)=\{x_1,...,x_n\},\quad  f^{-1}(ix)=\{y_1,...,y_n\}.
$$ If we denote for $i,j = 1, \ldots , n$
$$
P_{ij} = x_i + y_j \in X \subset \tilde{C}^{(2)}
$$
then  $\pi^{-1}(z)=\{P_{ij}\,; \,\, i,j = 1, \ldots, n\}$.\\
By construction the image $D(P_{ij}) = (p_2)_*(D \cap (\{P_{ij}\} \times X))$ is given by

$$D(P_{ij})=\sum_{l=1, l \neq j}^n P_{il} + \sum_{k=1, k \neq i}^n P_{kj}$$

In particular the correspondence $D$ is fixed point free. Moreover a straightforward computation shows
that

$$D^2(P_{ij})=(2n-2)P_{ij} + (n-2)D(P_{ij}) + 2 \sum_{k\neq i, l\neq j}P_{kl},$$

\noindent thus we get:

$$D^2(P_{ij}) - (2n-4)P_{ij} - (n-4)D(P_{ij})=2 \pi^*(\pi(P_{ij})).$$

\noindent This implies that the endomorphism $\gamma_D$ of the Jacobian $J(X)$ induced by $D$ satisfies the equation
$$
\gamma_D^2 + (4-n)\gamma_D - (2n-4)=0.
$$

Recall that according to a theorem of Kanev (see \cite{lb}, Theorem 12.9.1) an effective fixed point free symmetric correspondence $D$ on a smooth projective curve
$X$ defines a Prym-Tyurin variety of exponent $e$ if and only if the endomorphism $\gamma_D$ associated to $D$ satisfies the equation
$$
\gamma_D^2 + (e-2)\gamma_D - (e-1)=0.
$$
Together with the above reasoning this implies:

\begin{prop}
$P = Im (\gamma_D - 1)$ is a Prym-Tyurin variety for the curve $X$ if and only if $n = 3$. In this case the exponent of $X$ is 3.   
\end{prop}

\begin{rem}
Observe that on $X$ there is another natural effective symmetric correspondence:

$$D'(P_{ij})=\sum_{k\neq i, l\neq j}P_{kl}$$

\noindent which is symmetric and whose associated endomorphism $\gamma_D \in \End(JX)$ satisfies the equation
$\gamma_{D'}^2 + (n-2)\gamma_{D'} -
(n-1)=0$. However $D'$ has fixed points, so we do not know whether it induces a Prym-Tyurin variety.
\end{rem}

\section{Existence of the curve $X$}

Let $h: C \lra \P^1$ be a hyperelliptic covering of genus $g$ as above. We want to determine those \'etale coverings
$f: \tilde{C} \lra C$ of degree 3 for which the curve $X$ defined by diagram (2.1) is smooth and irreducible.\\
 
Recall that if $B_h = \{a_1, \ldots, a_{2g+2}\} \subset \P^1$ denotes the branch locus of $h$ and $\sigma_i$ denotes the class 
of the path from a fixed point $z_0 \in \P^1$ going around $a_i$ once, then
$$
\pi_1(\P^1 \setminus B_h, z_0) = <\sigma_1, \ldots, \sigma_{2g+2} : \prod^{2g+2}_{i=1} \sigma_i = 1 >.  \eqno(3.1)
$$ 
Let 
$$
\mu: \pi_1(\P^1 \setminus B_h, z_0) \lra S_6
$$
 be a classifying homomorphism for the composed map 
$f \circ h : \tilde{C} \lra C \lra \P^1$ and denote
$$
G = Im(\mu) \subset S_6.
$$ By construction $\mu(\sigma_i) = t_1t_2t_3$ where $t_1, t_2$ and $t_3$ are disjoint transpositions, 
but not all such products can occur.\\
In fact, if we denote as above $h^{-1}(z) = x + ix, \,\,\, f^{-1}(x) = \{x_1, x_2,x_3\}$ and $f^{-1}(ix) = \{y_1,y_2,y_3\}$ and if we 
identify $(x_1,x_2,x_3)$ with $(1,3,5)$ and $(y_1,y_2,y_3)$ with $(2,4,6)$, then exactly the following 6 permutations can occur 
$$
\begin{array}{c} \{(1\; 2)(3\; 4)(5\; 6),
(1\; 4)(2\; 5)(3\; 6),(1\; 6)(2\; 3)(4\; 5),\\
(1\; 2)(3\; 6)(4\; 5), (1\; 4)(2\; 3)(5\; 6),(1\; 6)(2\; 5)(3\; 4).\}
\end{array}
$$
Hence the existence of an \'etale covering $f: \tilde{C} \lra C$ is equivalent 
to the existence of a homomorphism $\mu: \pi_1(\P^1 \setminus B_h, z_0) \lra S_6$ as above such that the image $G = Im(\mu)$ is a 
transitive subgroup of $S_6$. By a direct computation one checks that there are up to conjugation exactly 3 types 
of transitive subgroups of $S_6$ generated by a subset of the above set of permutations, namely:\\

\noindent
I: $G = \;<(1\; 2)(3\; 4)(5\; 6), (1\; 4)(2\; 5)(3\; 6)> \,\, \simeq S_3,$\\
                            
\noindent
II: $G = \;< (1\; 2)(3\; 4)(5\; 6), (1\; 4)(2\; 5)(3\; 6), (1\; 2)(3\; 6)(4\; 5)>\,\, \simeq S_2 \times S_3$\\

\noindent  
III: $G$ generated by all 6 permutations of above, i.e.\\
\hspace*{0.6cm}  $G = \;< (2 \; 4 \;6), (1 \; 5)(2 \;4),  (1\; 4)(2\; 5)(3\; 6)>\; \simeq S_3 \times S_3$.\\

In order to see in which cases the associated curve $X$ is smooth and irreducible, we describe the monodromy 
associated to the constuction of $X$. For this we have to analyze the action of the group $G$ on a fibre of $\pi: X \lra \P^1$,
i.e. the action of $G$ on the set $\{x_1 +y_1, x_1 +y_2, \ldots, x_3+y_3\}$. One immediately checks that the in the cases I and II 
this action is not transitive. Hence the normalization of $X$ is not connected in these cases.
So for the  rest of this section let $G$ denote the group of case III. Then we have:

\begin{lem}
If $f: \tilde{C} \lra C$ is an \'etale covering of degree 3 of a hyperelliptic curve $C$ such that the image of a classifying homomorphism
$\mu: \pi_1(\P^1 \setminus B_h, z_0) \lra S_6$ is a group of type III, then the curve $X$ of diagram (2.1) is smooth and irreducible.
\end{lem}

\begin{pf} The stabilizer of the element $P_{11} = x_1 + y_1$ of the fibre $\pi^{-1}(z)$ is the group $G_{P_{11}} = < (1\; 2)(3\; 4)(5\; 6),
(3\; 4)(5\; 6) >$, which is Klein's group of 4 elements. Since $G$ is of order 36, this means that $G$ acts transitively on the set 
$\{x_1 +y_1, x_1 +y_2, \ldots, x_3+y_3\}$ implying that $X$ is irreducible. The proof of the fact that $X$ is smooth is a slight 
generalization of the proof of \cite{lb}, Lemma 12.8.1.  $\hfill{\square}$
\end{pf}

\noindent
Let 
$$ 
\xymatrix@R=0.4cm@C=0.4cm{
\mbox{Gal}(h \circ f) \ar[ddd]_{\gamma} \ar[dr]\\
& \tilde{C} \ar[d]^f\\
& C \ar[dl]^h\\
\P^1          }
$$
denote the Galois extension of $h \circ f: \tilde{C} \lra \P^1$. 
The next proposition identifies $\mbox{Gal}(h \circ f)$ in terms of the geometric construction of section 2. 
For this consider the curve $Y = \tilde{D} - \tilde{\Delta} \subset \tilde{X} \times \tilde{X}$ of section 1. It can be considered 
as a symmetric (2,2)-correspondence on $\tilde{X}$ without fixed points. So the projections $q_1$ and $q_2:Y \lra \tilde{X}$
coincide and are \'etale of degree 2.  Let us observe that in this case $(\deg f =3)$ the map $p_1: \tilde{X} \lra \tilde{C}$
is \'etale (see proof of Lemma 5.1). Let $\delta$ denote the composed map $Y \stackrel{q_1}{\lra} \tilde{X} \lra \P^1$. 

\begin{prop}
The map $\delta:Y \lra \P_1$ coincides with the Galois extension $\gamma: \mbox{Gal}(h \circ f) \lra \P_1$.
\end{prop}

\begin{pf}
As in section 2, for $z \in \P^1$, denote $h^{-1}(z)=x+ix$ and $f^{-1}(x)=\{x_1,x_2,x_3\},\quad  f^{-1}(ix)=\{y_1,y_2,y_3\}.$
Then the fibre $\delta^{-1}(z)$ for a general $z \in \P^1$ consists of the 36 elements 
$\{((x_1, y_1), (x_1, y_2)), ((x_1, y_1), (x_1, y_3)), \newline
((y_1,x_1), (y_1, x_2)), ((y_1,x_1), (y_1, x_3)), \ldots , ((y_3,x_3), (y_3, x_2)) \}$ It is immediate to check that the group 
$G$ acts transitively on these fibres or equivalently that the stabilizer of a point, say $((x_1, y_1), (x_1, y_2))$, is trivial. 
This implies that the Galois covering $\mbox{Gal}(h \circ f)$ is a normalization of $Y$. The smoothness of $Y$ follows from the 
fact that $Y$ is a symmetric fixed point free correspondence on the smooth curve $\tilde{X}$. $\hfill{\square}$
\end{pf}

In order to study the Prym-Tyurin variety $P$ of Proposition 2.1 we have to take into account also the subgroups of $G$, since
to every such subgroup there corresponds an intermediate covering of $\delta:Y \lra \P_1$. 
For this note that the triple products of transpositions generating the group $G$ form two conjugation classes in $G$, namely
$$
C_1 = \{(1\; 2)(3\; 4)(5\; 6),
(1\; 4)(2\; 5)(3\; 6),(1\; 6)(2\; 3)(4\; 5)\}. 
$$
$$
C_2 = \{(1\; 2)(3\; 6)(4\; 5), (1\; 4)(2\; 3)(5\; 6),(1\; 6)(2\; 5)(3\; 4)\},
$$
Moreover we have 
$$
G = \; <C_1,C_2> \; \simeq \; <C_1> \times <C_2> \; \simeq S_3 \times S_3.  \eqno(3.2)
$$
Let $\tau: G \lra G$ denote the outer automorphism interchanging the direct factors $<C_1>$ and $<C_2>$ of $G$. The subgroup 
diagram of $G$ consists of 8 conjugacy classes of subgroups invariant under $\tau$ and of 14 pairs of different conjugacy
classes of subgroups $(G_1,G_2)$ with $\tau(G_1) = G_2$. We need only the following part of it

$$
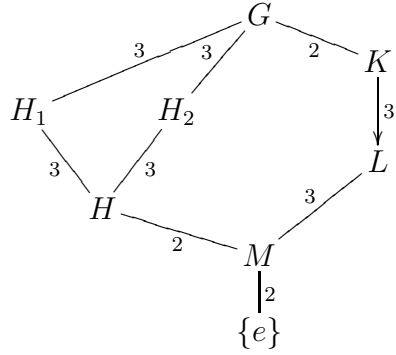
\xymatrix@R=0.2cm@C=0.4cm{
     & & & G \ar@{-}[ddlll]_3 \ar@{-}[ddl]_3 \ar@{-}[drr]_2\\
     & & & & &K \ar[dd]^3\\
   H_1 \ar@{-}[ddr]_3 & & H_2 \ar@{-}[ddl]^3 \\
     & & & & &L \ar@{-}[ddll]_3\\
     & H \ar@{-}[drr]_2\\
     & & & M \ar@{-}[dd]^2\\ 
     & & & & &\\     
     & & & \{e\}   }
$$

\noindent
where \\

$\begin{array}{ll}
      H &= \; < (1\; 2)(3\;4)(5\,6),(3\;5)(4\;6) >,\\
      H_1 &= \; < (1\; 2)(3\;4)(5\,6),(3\;5)(4\;6), (1\;3\;5)(2\;4\;6)>\\
      &= \; < (1\; 2)(3\;4)(5\,6), C_2>,\\
      H_2 &= \; < (1\; 2)(3\;4)(5\,6),(3\;5)(4\;6), (1\;3\;5)(2\;6\;4)>\\
      &= \; <(1\; 2)(3\;6)(4\,5), C_1>,\\
      K &= \; < (2\;4\;6),(1\;3\;5),(1\;5)(2\;4) >,\\
      L &= \; < (2\;4\;6), (3\;5)(4\;6) > \; \mbox{and}\\
      M &= \; < (3\;5)(4\;6) >.
      \end{array}
      $\\
      
\noindent

\begin{lem}
Up to conjugation in $G$ we have\\
\mbox{(a)}: $Y/L = \tilde{C}$ and $Y/K = C$,\\
\mbox{(b)}: $Y/H = X$ and $Y/M = \tilde{X}.$
\end{lem}

\begin{pf} Observe first that $L = G \cap S_5$, where
$S_5$ denotes the stabilizer of the symbol 1 in $S_6$. Hence $\mu^{-1}(L) \subset \pi_1(\P^1 \setminus B_h, z_0)$
is isomorphic to the fundamental group of $\tilde{C} \setminus (h \circ f)^{-1}B_h$. Hence $Y/L = \tilde{C}$ and by 
Galois theory it follows that $Y/K = C$, since there is no other subgroup between $L$ and $G$. This completes the proof of (a).\\  
In order to see (b), note first that the subgroup $H$ is the stabilizer in $G$ of the set $\{1,2\}$. 
This implies that the action of $G$ on the set of classes
$G/M$ gives a homomorphism $\lambda: G \lra S_9$, which is injective and with image a transitive subgroup of $S_9$. 
Hence, by construction of the curve $X$ the composition 
$\lambda \circ \mu: \pi_1(\P^1 \setminus B_h, z_0) \lra S_9$ 
is a classifying morphism for the covering $\pi: X \lra \P^1$. Denoting by $S_8$ the stabilizer of a symbol in $S_9$, this implies that 
$(\lambda \circ \mu)^{-1}(S_8)$ is isomorphic to the fundamental group of $X \setminus \pi^{-1}(B_h)$. But a short computation shows that 
$\lambda^{-1}(S_8) = H$. Hence the fundamental group of $X \setminus \pi^{-1}(B_h)$ is isomorphic to $\mu^{-1}(H)$ implying 
$Y/H = X$. Since $M$ is the the stabilizer of the ordered set $(1\,,2)$ in $G$, one shows in a similar way that 
$Y/M = \tilde{X}.$ $\hfill{\square}$
\end{pf}

Combining everything and denoting $X_1 = Y/H_1$ and $X_2 = Y/H_2$, we obtain the following diagram of morphisms of 
smooth projective curves:
             
$$
\xymatrix@R=0.1cm@C=0.4cm{
    & & &  Y \ar[dd]^{2:1}\\
    & & & & & \\
    & & &  \tilde{X} \ar[dll]_s^{2:1} \ar[ddrr]\\
    & X \ar[ddl]^{3:1} \ar[ddddrr]_{\pi} \ar[ddrr]\\
    & & & & &\tilde{C} \ar[dd]_f^{3:1}\\
    X_1 \ar[ddrrr]_{3:1}^{f_1} &&& X_2 \ar[dd]^{f_2}\\
    & & & & & C \ar[dll]_h^{2:1}\\
    & & & \P^1    }    \eqno(3.3)
$$\\    

\noindent
Recall from section 2 the correspondence $D = (\lambda \times \lambda)_*(Y) \subset X \times X$. 

\begin{prop}
$\lambda \times \lambda|_Y:Y \lra D$ is an isomorphism.
\end{prop}

\begin{pf} The map is given by $((x_i,y_j),(x_i,y_k)) \lra (x_i+y_j,x_i+y_k)$. The morphism $D \lra Y$ defined by 
$(a,b) \lra ((a\cap b, a-a\cap b), (a\cap b, b-a\cap b))$ is inverse to it. Here $a$ and $b$ are considered as divisors on 
the curve $\tilde{C}$ and $a \cap b$ denotes the greatest common divisor of $a$ and $b$.  $\hfill{\square}$
\end{pf}

\section{The moduli spaces}

Let the notation be as at the end of the last section. So $f: \tilde{C} \lra C$ is an \'etale covering of degree 3 
of the hyperelliptic curve $h: C \lra \P^1$ with branch locus $B_h = \{a_1, \ldots, a_{2g+2}\}$, 
such that the composition $h \circ f$ is given by a classifying homomorphism 
$\mu: \pi_1(\P^1 \setminus B_h, z_0) \lra G \subset S_6$ as in (3.1) and 
$G \simeq S_3 \times S_3$, of type III. In this section we want to study the moduli of this situation.\\

Consider again the direct product decomposition (3.2).
Since $\mu(\sigma_i) \in C_1 \cup C_2$, we can enumerate the $\sigma_i$ in such a way that 
$$
\mu(\sigma_i) = (g_i,1) \in \;<C_1> \times <C_2> \quad \mbox{for} \quad i=1, \ldots, \alpha,
$$
$$
\mu(\sigma_i) = (1,g_{\alpha+i}) \in \;<C_1> \times <C_2> \quad \mbox{for} \quad i=1, \ldots, \beta
$$
for some $\alpha$ and $\beta$ with $\alpha + \beta = 2g+2$. In particular $g_i \in C_1$ for $i=1,\ldots, \alpha$ and 
$g_i \in C_2$ for $i=\alpha + 1, \ldots, \alpha + \beta$. Then the condition $\prod^{2g+2}_{i=1} \sigma_i = 1$ 
is equivalent to the two conditions
$$
\prod_{i=1}^{\alpha} g_i =1 \quad \mbox{and} \quad \prod_{i=1}^{\beta} g_{\alpha+i} = 1   \eqno(4.1)
$$  
So we must have $\alpha$ and $\beta$ even with
$$\alpha \geq 4 \quad \mbox{ and} \quad  \beta \geq 4,
$$ 
since under the isomorphisms $<C_1>\; \simeq S_3$ and $<C_2>\; \simeq S_3$ 
the elements of $C_1$ and $C_2$ correspond to transpositions, that is are of order 2.\\ 
Using the notation of above it makes sense to call an \'etale covering $f: \tilde{C} \lra C$ of degree 3 
{\it of type $(\alpha, \beta)$} if $\mu(\sigma_i) \in C_1$ for $i=1, \ldots, \alpha$ and $\mu(\sigma_i) \in C_2$ for 
$i = \alpha + 1, \ldots, \alpha + \beta= 2g+2$. The following theorem is
the main result for studying the moduli space of \'etale degree 3 coverings of hyperelliptic curves of type $(\alpha, \beta)$.  

\begin{thm}
Suppose $\alpha, \beta$ are even integers $\geq 4$ with $\alpha + \beta = 2g+2$ and $a_1, \ldots, a_{2g+2} \in \P^1$ 
pairwise different. There is a canonical $1:1$-correspondence
between the sets of\\
\begin{enumerate}
\item coverings $\tilde{C} \stackrel{f}{\lra} C \stackrel{h}{\lra} \P^1 \; \; \deg h =2$ ramified exactly over $a_1, \ldots, a_{2g+2}$ and 
      $f$ unramified of type $(\alpha, \beta)$ of degree 3,  and
\item pairs of trigonal curves $f_1: X_1 \lra \P^1$ simply ramified exactly over $a_1, \ldots, a_{\alpha}$ and $f_2: X_2 \lra \P^1$ 
simply ramified exactly over $a_{\alpha+1}, \ldots, a_{2g+2}$.
\end{enumerate}         
\end{thm}

\begin{pf}
Given $(f,h)$ as in (1) it follows from (4.1) that the homomorphism $\mu$ induces homomorphisms $\mu_1$ and $\mu_2$ such that the 
following diagram is commutative
$$
\xymatrix@R=0.6cm@C=1cm{
     \pi_1(\P^1 \setminus \{a_1, \ldots, a_{\alpha}\},z_0) \ar[r]^-{\mu_1} & S_3\\
     \pi_1(\P^1 \setminus \{a_1, \ldots, a_{2g+2}\},z_0)  \ar[u]^{l_1} \ar[r]^-{\mu} \ar[d]_{l_2} &
     G = \;<C_1> \times <C_2> \ar[u]_{p_1}  \ar[d]^{p_2}\\
     \pi_1(\P^1 \setminus \{a_{\alpha+1}, \ldots, a_{2g+2}\},z_0) \ar[r]^-{\mu_2} & S_3     }  \eqno(4.2)
$$  
where $l_1$ is induced by the inclusion $\P^1 \setminus B_h \lra \P^1 \setminus \{a_1, \ldots , a_{\alpha}\}$ and similarly
$l_2$ by $\P^1 \setminus B_h \lra \P^1 \setminus \{a_{\alpha+1}, \ldots , a_{2g+2}\}$ and $p_i: G \lra <C_i> \; \simeq S_3$ 
for $i=1$ and 2 is the projection followed by a fixed isomorphism $<C_i> \; \simeq S_3$.\\
To $\mu_1$ corresponds a triple covering $f_1:X_1 \lra \P^1$ simply ramified exactly over $a_1, \ldots, a_{\alpha}$. Similarly 
to $\mu_2$ corresponds a triple covering $f_2:X_2 \lra \P^1$ simply ramified exactly over $a_{\alpha+1}, \ldots, a_{2g+2}$. 
The coverings $f_1$ and $f_2$ are uniquely determined up to an automorphism, since every automorphism of $S_3$ is inner.\\
Conversely, given $(f_1,f_2)$ as in (2), to $f_1$ corresponds a homomorphism $\mu_1$ and to $f_2$ a homomorphism $\mu_2$ as in the 
diagram. Define $i_1: S_3 \lra G$ as the composition of a fixed isomorphism $S_3 \simeq \; <C_1>$ with the embedding of $<C_1>$ 
as the first factor and similarly $i_2: S_3 \lra G$ as the composition of a fixed isomorphism $S_3 \simeq \; <C_2>$ with 
the embedding of $<C_2>$ as the second factor. Then we can define $\mu: \pi_1(\P^1 \setminus \{a_1, \ldots, a_{2g+2}\},z_0) \lra G$
such that $\mu = i_1 \circ \mu_1 \circ l_1 = i_2 \circ \mu_2 \circ l_2$ by setting
$$
\mu(\sigma_i) = i_1 \circ \mu_1 \circ l_1(\sigma_i) \quad \mbox{for} \quad i= 1, \ldots \alpha \quad \mbox{and}
$$
$$
\mu(\sigma_i) = i_2 \circ \mu_2 \circ l_2(\sigma_i) \quad \mbox{for} \quad i= \alpha +1, \ldots 2g+2.
$$
Since $i_1$ and $i_2$ map the transpositions of $S_3$ to the products of transpositions in $<C_1>$ and $<C_2>$ i.e. 
to $C_1$ and $C_2$, the homomorphism $\mu$ defines a composition of coverings $\tilde{C} \stackrel{f}{\lra} C \stackrel{h}{\lra} \P^1$
as in (1). Certainly these maps are inverse to each other which completes the proof.  $\hfill{\square}$  \end{pf}

\begin{prop}
The maps $f_i: X_i \lra \P^1$ of theorem 4.1 and of diagram (3.3) coincide. In particular
$$X_i = Y/H_i \quad \mbox{ for} \quad i= 1, 2.$$ 
\end{prop}

\begin{pf}
Observe that if we fix the isomorphism $S_3 \stackrel{\sim}{\rightarrow} \;<C_1>$ by
$$
(1\,2) \mapsto (1\,2)(3\,4)(5\,6) \quad (1\,3) \mapsto (1\,4)(2\,5)(3\,6),
$$
and $S_3 \stackrel{\sim}{\rightarrow} \;<C_2>$ by
$$
(1\,2) \mapsto (1\,2)(3\,6)(4\,5), \quad (1\,3) \mapsto (1\,4)(2\,3)(5\,6)
$$
then with the notation of diagram (4.2) we have that the group $\mu^{-1}(\{1,(1\,2)\})$ is isomorphic to the 
fundamental group of $X_1 \setminus f_1^{-1}(\{a_1, \ldots,a_{\alpha}\})$. Hence the group $l^{-1}_1\mu_1^{-1}(\{1,(1\, 2)\})$ 
is isomorphic to the fundamental group of $X_1 \setminus f_1^{-1}(\{a_1, \ldots, a_{2g+2}\})$. 
On the other hand 
$$
l_1^{-1}(\mu_1^{-1}(\{1,(1\,2)\}) = \mu^{-1}(<(1\,2)(3\,4)(5\,6)> \times <C_2>) = \mu^{-1}(H_1)
$$
which gives $X_1 \simeq Y/H_1$. Similarly 
$$
l_2^{-1}(\mu_2^{-1}(\{1,(1\,2)\}) = \mu^{-1}(<C_1> \times  <(1\,2)(3\,6)(4\,5)>) = \mu^{-1}(H_2)
$$
which gives $X_2 \simeq Y/H_2$.  $\hfill{\square}$
\end{pf}

\begin{cor}
Let $\alpha, \beta$ be even integers $\geq 4$ with $\alpha + \beta = 2g+2$ and $g \geq 3$.\\
The moduli space $\M(\alpha, \beta)$ of \'etale degree 3 coverings of hyperelliptic curves of genus $g$ of type $(\alpha,\beta)$ is isomorphic 
to the moduli space of pairs of trigonal coverings of genus $\frac{\alpha}{2} - 2$ and $\frac{\beta}{2}-2$ with 
simple ramification and disjoint branching. In particular $\dim \M(\alpha,\beta) = \alpha+\beta-3 = 2g-1$.  
\end{cor}

\begin{pf}
Theorem 4.1 or to be more precise a slight generalization of it concerning families of the corresponding coverings implies that the 
moduli functors in question are isomorphic. This implies the statement about the moduli spaces. It remains to compute the dimension.
Consider 
$$
A := (\times_{i=1}^{\alpha} \P^1) \setminus \Delta_{\alpha} \qquad \mbox{and} \qquad B:= (\times_{i=1}^{\beta} \P^1) \setminus \Delta_{\beta}
$$
where $\Delta_{\alpha}$ and $\Delta_{\beta}$ denote the corresponding discriminants and let 
$$
\pi_{\alpha}: H^{3,\alpha} \lra A \quad \mbox{and} \quad \pi_{\beta}: H^{3,\beta} \lra B
$$
 denote the Hurwitz spaces of triple covers of $\P^1$ simply ramified in $\alpha$ respectively $\beta$ points.
Let $R \subset A \times B$ be the open set
$$
R = \{(a_1,\ldots,a_{\alpha},a_{\alpha+1},\ldots, a_{\alpha+\beta}) \;|\; a_i \not= a_{\alpha+j} \; \mbox{for} \; 1 \leq i \leq \alpha, 
1 \leq j \leq \beta \}
$$ 
Then $(\pi_{\alpha} \times \pi_{\beta})^{-1}(R) \lra R$ parametrizes pairs of simple triple covers with disjoint branching. 
The action of $PGL(1)$ on $R$ lifts to an action on $(\pi_{\alpha} \times \pi_{\beta})^{-1}(R)$ and the quotient
$(\pi_{\alpha} \times \pi_{\beta})^{-1}(R)/PGL(1)$ is the moduli space of pairs of trigonal coverings with simple 
ramification and disjoint branching. So $\dim \M(\alpha,\beta) = \alpha + \beta - 3$. $\hfill{\square}$
\end{pf}

\section{Comparison with other ppav's}

It is the aim of this section to relate the Prym-Tyurin varieties $P$ introduced in section 2 to the other 
principally polarized abelian varieties occurring in this situation. For this we first 
determine the ramification of the maps and the genera of the curves occurring in diagram (3.3).\\

\noindent
 Recall that $C$ is a hyperelliptic curve of genus $g$ 
ramified over the points $a_1, \ldots, a_{2g+2} \in \P^1$ and $\tilde{C}$ an \'etale covering of degree 3 of $C$ and thus 
of genus $3g-2$. Moreover from the construction we get that 
all maps in the diagram (3.3) are unramified over $\P^1 \setminus\{a_1, \ldots , a_{2g+2}\}$. Suppose $f:\tilde{C} \lra C$
is of type $(\alpha,\beta)$. Then we have

\begin{lem}
(a) $g(\tilde{C}) = 3g-2, \; g(\tilde{X}) = 9g-8, \; g(Y) = 18g-17;$\\
(b) $ g(X) = 3g-5$;\\
(c) $g(X_1) = \frac{\alpha}{2} - 2$
and $g(X_2) = \frac{\beta}{2} - 2$; \\
(d) $\dim P = g-3$.
\end{lem} 

\begin{pf}
The covering  $\tilde{X} \lra C$ is defined by restricting the \'etale map $f \times f$ in the following diagram
$$ 
\xymatrix@R=0.4cm@C=0.4cm{
\tilde{X} = (f \times f)^{-1}(C) \ar[d] \ar @{^{(}->} [r] & \tilde{C}^2 \ar[d]^{f \times f}\\
C = \{(x,ix)\;|\;x \in C\} \ar @{^{(}->} [r] & C^2  }
$$
It follows that $\tilde{X} \lra \tilde{C}$ is \'etale of degree 3. Recall from section 2 that
$Y$ is a fixed point free symmetric $(2,2)$-correspondence on the curve $\tilde{X}$. Hence $Y \lra \tilde{X}$ is an \'etale 
double covering and Hurwitz formula gives (a). From the description of the fibre $(fh)^{-1}(z)$ for any $z \in \P^1$ 
in section 2 we see that $\pi$ is of ramification type $(2,2,2,1,1,1)$. Again Hurwitz formula gives (b). 
(c) was proven already in Corollary 4.3.\\
Proof of (d): Let $N_P$ denote the norm endomorphism associated to the endomorphism $\gamma_D$ of the correspondence $D$.
According to \cite{lb}, 5.3.10 $\dim P$ is related to the analytic trace of $N_P$ by
$$
\dim P = \frac{1}{3} \Tr_a(N_P).
$$
Since $\gamma_D = 1_X - N_P$, we have 
$$
\Tr_r(\gamma_D) = 2g(X) - \Tr_r(N_P).
$$
where $\Tr_r$ denotes the rational trace, which is related to the analytic trace by $\Tr_r = 2 Re \Tr_a$.
On the other hand, according to a theorem of Weil (see \cite{lb}, 11.5.2 and 3.1.3) we have 
$$
\Tr_r(\gamma_D) = 8,
$$
since the correspondence $D$ is without fixed points. Putting everything together gives the assertion.  $\hfill{\square}$
\end{pf}

\noindent
From Lemma 5.1 we deduce 
$$
\dim f_1^*JX_1 + \dim f_2^*JX_2 = \dim P.
$$
This suggests that there should be a relation between $JX_1 \times JX_2$ and $P$. The following theorem 
is the main result of this section.

\begin{thm}
The canonical map $f_1^* + f_2^* : JX_1 \times JX_2 \lra JX$ induces an isomorphism
$$
JX_1 \times JX_2 \simeq P.
$$
\end{thm}

\noindent
The proof consists of a careful analysis of the action of the group $G$. Recall that $S_3$ admits 3 (absolutely) 
irreducible $\Q$-representations, the trivial and the alternating representations $1_{S_3}$ and $U$ of dimension 1 
and the standard 2-dimensional representation $V$. The tensor products of these are all irreducible representations of 
the group $G = S_3 \times S_3$. For any subgroup $H$ of $G$ let $\rho_H(1_H)$ denote the representation of $G$ 
induced by the trivial representation $1_H$ of $H$.

\begin{lem}
\hspace*{1.3cm} $\rho_{H_1}(1_{H_1}) - 1_G \simeq V \otimes 1_{S_3}$,
$$
\rho_{H_2}(1_{H_2}) - 1_G \simeq 1_{S_3} \otimes V.
$$
\end{lem}

\begin{pf}
Fix isomorphisms $<C_1> \; \simeq S_3$ and $<C_2> \; \simeq S_3$ and thus $G \simeq S_3 \times S_3$. Then $H_1 \simeq 
\; <(1\;2)(3\;4)(5\;6)> \times <C_2> \, \simeq S_2 \times S_3$.\\
The representation $\rho_{H_1}(1_{H_1})$ is then given by the action of the group $G = S_3 \times S_3$ on the quotient
$(S_3 \times S_3)/(S_2 \times S_3)$ or equivalently by the action of the first factor
$S_3$ on the quotient $S_3/S_2$. But it is easy to see that this is just the representation $V + 1_{S_3}$. This implies the 
first equation. The second equation is proved in the same way. $\hfill{\square}$
\end{pf}  

\begin{pf} ({\it of Theorem 5.2}) The action of the group $G$ on the curve $Y$ induces a homomorphism 
$\Q[G] \lra End_0(JY)$ of the rational group ring $\Q[G]$. Using this one can associate an abelian subvariety of $JX$ 
to every projector $p \in \Q[G]$ and thus to every subrepresentation of $\Q[G]$ in a natural way. 
In particular, if $W$ is an irreducible $\Q$-representation of $G$, and $( \;\;,\,\;)$ a $G$-invariant scalar product of $W$,
then for any nonzero $w \in W$ a projector for $W$ is given as follows (see \cite{lb}, p. 434)
$$
p_w = \frac{\dim W}{|G| \cdot \Vert w \Vert^2} \sum_{g \in G} (w,gw)g
$$
It is well known 
(see e.g. \cite{lr2}, Corollary 3.2), that the pull-back of the Prym variety $P(f_i)$ of the morphism 
$f_i: X_i = Y/H_i \lra Y/G = \P^1$ in $J_Y$ corresponds to the 
representation $\rho_{H_i}(1_{H_i}) - 1_G$.
But according to Lemma 5.3 both representations are irreducible. So choosing nonzero vectors 
$w_1 \in V \otimes 1_{S_3}$ and 
$w_2 \in 1_{S_3} \otimes V$ and denoting $\varphi$ the composed map $Y \lra \tilde{X} \stackrel{s}{\lra} X$ of diagram (3.3),
this implies
$$
p_{w_i}(JY) \sim \varphi^*f_i^*JX_i  \eqno(5.1)
$$
for $i=1$ and 2, where $\sim$ denotes isogeny. \\
Now an easy computation shows that 
the vector spaces $(V \otimes 1_{S_3})^H$ and $(1_{S_3} \otimes V)^H$ of $H$-invariants are one-dimensional. 
Choosing $w_1 \in (V \otimes 1_{S_3})^H \setminus 0$ and $w_2 \in (1_{S_3} \otimes V)^H \setminus 0$ we have
$$
(w_i,hgw_i) = (w_i,ghw_i) =(w_i,gw_i)
$$
for any $h \in H$ and $i=1,2$. Hence applying \cite{lr1}, Proposition 3.5 we get that the projector $p_{w_i}$ descends 
to an endomorphism 
$
\tilde{p}_{w_i} : JX \lra JX
$
such that 
$$
p_{w_i}(JY) = \varphi^*(\tilde{p}_{w_i}(JX)).
$$
In fact 
$$
\tilde{p}_{w_i}(x) = \frac{\dim V}{|G|\cdot \Vert w \Vert^2} \sum_{g \in G} (w,gw) \bar{g}(x),
$$
where 
$$\bar{g}(x) = \sum_{h \in H} \varphi(ghy)  \eqno(5.2)
$$
 with $y \in Y$ such that $\varphi(y)=x$.\\
So (5.1) implies
$$
\tilde{p}_{w_i}(JX) \sim f_i^*JX_i  \eqno(5.3)
$$
Now the projector $\tilde{p}_{w_i}$ does not depend on the choice of the vector $w_i$ the vector spaces $(V \otimes 1_{S_3})^H$
and $(1_{S_3} \otimes V)^H$ being one-dimensional. This implies that in (5.3) we actually have equality instead of only isogeny  
(see \cite{lb}, Proposition 13.6.4). So we conclude
$$
(f_1^* + f_2^*)(JX_1 \times JX_2) = (\tilde{p}_{w_1} + \tilde{p}_{w_2})(JX) \eqno(5.4)
$$
In order to complete the proof of the theorem we have to show that the image of $\tilde{p}_{w_1} + \tilde{p}_{w_2}$ is just $P$.
For this we have to compute the projectors $p_{w_1}$ and $p_{w_2}$ explicitly.\\ 

\noindent
Consider the decomposition
$$
G = H_1 \cup (2\,4\,6)H_1 \cup (1\,3\,5)H_1
$$
If we define $v_1 = H_1, \; v_2 = (2\,4\,6)H_1, \; v_3 = (1\,3\,5)H_1$ and $e_1 = v_1 - v_2, \; e_2 = v_2 - v_3$, 
then the action of $G$ on $G/H_1$ gives us the representation $V \otimes 1_{S_3} = e_1\Q \oplus e_2\Q \simeq \Q^2$.\\ 
Let us for example consider the action of $(3\,5)(4\,6) \in H$. We have 
$(3\,5)(4\,6)(H_1) = H_1, \; (3\,5)(4\,6)((2\,4\,6)H_1) = (1\;3\;5)H_1$ and hence $(3\,5)(4\,6)(e_1) = e_1 + e_2$ and 
$(3\,5)(4\,6)(e_2) = - e_2$.
This means in matrix form $(3\,5)(2\,4) = \left( \begin{array}{cc} 1 &0\\ 1 & -1 \end{array} \right): \Q^2 \lra \Q^2$.\\
The next step is to introduce a $G$-invariant scalar product on $\Q^2$. One checks directly that
$$
<(u_1,u_2),(v_1,v_2)> = 2u_1v_1 -u_1v_2 - u_2v_1 + 2u_2v_2
$$
is a $G$-invariant scalar product on $V \otimes 1_{S_3} = \Q^2$. We noted already above that $\dim (V \otimes 1_{S_3})^H = 1$. 
In fact, one observes that the matrix representation for the elements of $H$ is either the above example or the identity, implying
that $(2,1)^{t} \in \Q^2$ is a generator of $(V \otimes 1_{S_3})^H$. \\
Choosing $w_1 = (2,1)^t$ we are ready to compute the projector
$p_{w_1}$. For example, the coefficient of $g = (2\,4)(3\,5)$ in $p_{w_1}$ is:
$$
\begin{array}{ll}
\frac{\dim(V \otimes 1_{S_3})}{|G| \cdot \Vert w \Vert^2}(w,gw) & 
= \frac{2}{36 \cdot 6} <(2,1),\left( \begin{array}{cc} -1&1\\0&1 \end{array} \right) \left( \begin{array}{c} 2 \\ 1 \end{array} \right) >\\
& = \frac{2}{36 \cdot 6}(-4-2+1+2) = \frac{2}{36 \cdot 6}(-3)  
\end{array}
$$
Proceeding in this way with all elements of $G$ and the representation $V \otimes 1_{S_3}$ and similarly the representation
$1_{S_3} \otimes V$ we obtain:
$$
\begin{array}{ll}
p_{w_1} + p_{w_2} \; = & \frac{2}{36 \cdot 6} \cdot \{12 \sum h  \\
& - \; 6 \sum (1\,3\,5)h - 6 \sum (1\,5\,3)h \\
& - \; 6 \sum (2\,4\,6)h - 6 \sum (2\,6\,4)h \\
& + \; 3 \sum (1\,3\,5)(2\,4\,6)h + 3 \sum (1\,3\,5)(2\,6\,4)h \\
& + \; 3 \sum (1\,5\,3)(2\,4\,6)h + 3 \sum (1\,5\,3)(2\,6\,4)h \}
\end{array}
$$
where the sum is always to be taken over all $h \in H$. (Notice that even though $H$ is not a normal subgroup of $G$,
the same expression for $p_{w_1} + p_{w_2}$ is valid if we use right cosets instead of left cosets.) So we get
$$
\begin{array}{ll}
36(\tilde{p}_{w_1} + \tilde{p}_{w_2}) \; =& 4\cdot 1_{JX} - 2 \overline{(1\,3\,5)} - 2 \overline{(1\,5\,3)} 
             - 2 \overline{(2\,4\,6)} - 2 \overline{(2\,6\,4)}\\
             & + \; \overline{(1\,3\,5)(2\,4\,6)} + \; \overline{(1\,3\,5)(2\,6\,4)}\\
             & + \; \overline{(1\,5\,3)(2\,4\,6)} + \; \overline{(1\,5\,3)(2\,6\,4)}
\end{array}
$$
where $\overline{(1\,3\,5)}$ is the endomorphism of $JX$ induced by $(1\,3\,5)$ using (5.2).\\
Now with the notation of identifications of sections 2 and 3 we have at the level of the curve $X$:
$$
\overline{(1\,3\,5)}(P_{11}) = P_{21}, \qquad \overline{(1\,5\,3)}(P_{11}) = P_{31},
$$
$$ 
\overline{(2\,4\,6)}(P_{11}) = P_{12}, \qquad \overline{(2\,6\,4)}(P_{11}) = P_{13},
$$ 
$$
\overline{(1\,3\,5)(2\,4\,6)}(P_{11}) = P_{22}, \qquad \overline{(1\,3\,5)(2\,6\,4)}(P_{11}) = P_{23},
$$
$$
\overline{(1\,5\,3)(2\,4\,6)}(P_{11}) = P_{32}, \qquad \overline{(1\,5\,3)(2\,6\,4)}(P_{11}) = P_{33}.
$$
So
$$
36(\tilde{p}_{w_1} + \tilde{p}_{w_2})(P_{11}) = 4P_{11} -2P_{21} -2P_{31} -2P_{12} -2P_{13} 
+ P_{22} + P_{23} + P_{32} + P_{33}
$$
On the other hand, recalling that $D(P_{11}) = P_{12} + P_{13} + P_{21} + P_{31}$ we get
$$
3(1_X - D)(P_{11}) + \pi^*(\pi(P_{11})) = 36(\tilde{p}_{V \otimes 1_{S_3}} + \tilde{p}_{1_{S_3} \otimes V})(P_{11})
$$
This implies
$$
(\tilde{p}_{w_1} + \tilde{p}_{w_2})(JX) = (1_{JX} - \gamma_D)(JX) = P
$$
which with (5.4) implies
$$
(f_1^* + f_2^*)(JX_1 \times JX_2) = P
$$
But on the one hand the restriction of the canonical principal polarization $\Theta$ of $JX$ to $P$ is of type (3,\ldots,3) 
and on the other hand the pull-back of $\Theta$ via $f_1^* + f_2^*: JX_1 \times JX_2 \lra JX $ is also of type (3,\ldots,3).
This implies that $f_1^* + f_2^*: JX_1 \times JX_2 \lra JX $ is a closed embedding, thus completing the proof of the theorem.$\hfill{\square}$
\end{pf}

\begin{rem}
We computed the dimension of all Jacobians and Prym varieties arising from the subgroup graph of the group $G$. 
It turns out that they are different from $\dim P$. So none of them is isogenous to the Prym-Tyurin variety $P$.
The details will not be included here.  
\end{rem}

\noindent
Now we are in a position to prove the Corollary stated in the introduction.\\

\noindent
{\it Proof} (of the Corollary of Theorems 4.1 and 5.2): Let $X_1, X_2$ and $X$ be as stated in the Corollary. 
According to Theorem 4.1 the pair of trigonal covers determines an \'etale degree 3 covering $f: \tilde{C} \lra C$ 
of a hyperelliptic curve $C$. So we are in the situation of Theorem 5.2 in the proof 
of which we saw that $f_1^* + f_2^*: JX_1 \times JX_2 \lra JX$ is an embedding. It only remains to be noted that the
curve $X$ of Theorem 5.2 coincides with the curve $X$ of the Corollary, i.e. is the fibre product of the trigonal covers. But this comes from the fact that
$$
<H_1,H_2> \;= G \qquad \mbox{and} \qquad H_1 \cap H_2 = H
$$ 
since the branchings of the trigonal covers are disjoint. $\hfill{\square}$

\section{The Abel-Prym map}

Let the notation be as above. In particular $X$ is a smooth projective curve of genus $3g-5$ and $D \subset X \times X$ the 
correspondence given by $D(P_{ij}) = \sum_{k \not= j} P_{ik} + \sum_{l \not=i}P_{lj}$ for $i,j = 1, \ldots, 3$. 
For any positive integer $d$ 
let $X^{(d)}$ denote the $d$-fold symmetric product of $X$ and $\alpha_d:X^{(d)} \lra JX$ be the Abel map with respect 
to a base point $p_0 \in X$. In this section we analyse the Abel-Prym map of $P$ which is by definition the composition
$$
\beta_P: X \stackrel{\alpha_1}{\lra} JX \stackrel{\gamma_D - 1_{JX}}{\lra} P \subset JX
$$   
By definition $\beta_P$ is given by the following diagram

$$
\xymatrix@C=0.4cm{
& X^{(4)} \times X \ar[dr]^{\alpha_4 - \alpha_1} \\
X \ar[ur]^{(D,1_X)} \ar[rr]_{\beta_P} && JX }
$$
where we consider the correspondence as a morphism $D: X \lra X^{(4)}$.

\begin{prop}
If $g\geq 6$, the Abel-Prym map $\beta_P: X \lra P$ is injective.
\end{prop}

\begin{pf} 
Suppose $\beta_P$ is not injective. Then according to the diagram there are two points $p_1, p_2 \in X$ such that
$D(p_1) - p_1 \sim D(p_2) - p_2$, where $\sim$ means linear equivalence. This implies
$$
D(p_1) + p_2 \sim D(p_2) + p_1,
$$
which means that $X$ admits a $g^1_5$, i.e. covering $X \lra \P^1$ of degree $\leq 5$. 
On the other hand $X$ is a three to one covering of a curve of genus $g_1 \leq \frac{g-3}{2}$, namely $X \lra X_1$.
But then Castenuovo's inequality (see \cite{cas}) implies
$$
3g-5 = g(X) \leq 2\cdot 4 + 3\cdot \frac{g-3}{2} \leq \frac{1}{2}(3g + 7).
$$
So $g \leq 5$, a contradiction.$\hfill{\square}$
\end{pf}

Let us now analyse the local behaviour of the Abel-Prym map $\beta_P$. The tangent space 
$T_0JX$ of $JX$ at 0 can and will be identified with $H^0(X,\omega_X)^*$. The differential $(d\gamma_D)_0:T_0JX \lra T_0JX$ has 
just 2 eigenvalues, namely 1 with multiplicity $g(X)-\dim P = 2g-2$ and $-2$ with multiplicity $\dim P = g-3$. Let $V_+$ and $V_-$
denote the corresponding eigenspaces. Clearly
$$
T_0P = V_-
$$
Define the {\it Prym-Tyurin canonical map} $\varphi_P$ of $P$ to be the composition of the canonical map 
$\varphi_X: X \lra \P(T_0JX)$ and the linear projection $r: \P(T_0JX) \lra \P(T_0P)$ with center $\P(V_+)$:
$$
\varphi_P: X \stackrel{\varphi_X}{\lra} \P(T_0JX) \setminus \P(V_+) \stackrel{r}{\lra} \P(T_0P).
$$ 
The following lemma is an immediate consequence of the fact that the canonical map $\varphi_X$ is everywhere defined:

\begin{lem}
A point $x \in X$ is a base point of the linear system $|T_0P| = |V_-|$ if and only if $\varphi_X(x) \in |V_+|.$
\end{lem}

For any point $p \in P$ there is a canonical isomorphism of tangent spaces $T_pP \simeq T_0P$ via which 
we identify the two vector spaces. In particular we consider the differential of $\beta_P$ at $x$ as a map
$(d\beta_P)_x: T_xX \lra T_0P$. Varying $x$ in $X$ we obtain a homomorphism of the corresponding tangent bundles
$d\beta_P:  T_X \lra  T_P = P \times T_0P$. It is easy to see that the projectivization of $d\beta_P$ 
coincides with the Prym-Tyurin canonical map
$$
\varphi_P = P(d\beta_P): X \lra \P(T_0P).
$$ 
Using these facts we are in a position to show:

\begin{prop}
If $g \geq 5$, the differential $(d\beta_P)_x: T_xX \lra T_0P$ of the Abel-Prym map is injective for every $x \in X$. 
\end{prop}

\begin{pf}
Assume that for some $x \in X$ the differential of the Abel-Prym map at $x$ is not injective, i.e. $(d\beta_P)_x = 0$. 
According to Lemma 6.2 and by what we have said above this means that $x$ is a base point of the linear system 
$|T_0P|$ which is the case if and only if the image of the canonical map satisfies $\varphi_X(x) \in \P(V_+)$. 
We will show that this leads to a contradiction. For this consider the commutative diagram
$$
\xymatrix@C=0.4cm{
X \ar[r]_D \ar[d]_{\alpha_1} & X^{(4)} \ar[d]^{\alpha_4}\\
JX \ar[r]^{\gamma_D} & JX }
$$
On the level of tangent spaces this gives
$$
\xymatrix@C=0.4cm{
T_xX \ar[r]_{(dD)_x} \ar[d]_{(d\alpha_1)_x} & T_{D(x)}X^{(4)} \ar[d]^{(d\alpha_4)_{D(x)}}\\
T_0JX \ar[r]^{(d\gamma_D)_0} & T_0JX }
$$
In the same way as in the proof of Proposition 6.1 we conclude from Castelnuovo's inequality that 
for $g(X)\geq 10$ or equivalently $g \geq 5$ the curve $X$ does not admit a $g^1_4$.
This implies that the map $(d\alpha_4)_{D(x)}: T_{D(x)}X^{(4)} \lra T_0JX$ is an isomorphism onto its image.\\
Now let $t \in T_xX$ be a nonzero vector. Since the projectivized differential of the Abel map $\alpha_1$ is the canonical 
map $\varphi_X$, we get from the assumption on $x$ that $d\alpha_1(t) \in V_+$. But $(d\gamma_D)_0$ is the identity on $V_+$,
implying 
$$
(d\gamma_D)_0(d\alpha_1)_x(t) = (d\alpha_1)_x(t)
$$ 
The commutativity of the above diagram implies
$$
(d\alpha_1)_x(t) \in (d\alpha_4)_{D(x)}T_{D(x)}X^{(4)} \subset T_0JX
$$
Projectivizing and using the fact that the projectivization of $(d\alpha_1)_x$ is the canonical map we get
$$
\varphi_X(x) \in \overline{D(x)} \subset |\omega_X|^*.   \eqno(6.1)
$$
Here $\overline{D(x)}$ means the linear span of the divisor $D(x)$ in the projective space $|\omega_X|^* = \P^{3g-6}$.
Now $h^0(D) = 1$ since $X$ does not admit a $g^1_4$. Hence the geometric version of Riemann-Roch (see \cite{acgh}, p. 12) implies
$$
\dim \overline{D(x)} = 3.
$$
But then (6.1) implies $\dim \overline{D(x) + x} = 3$ which again by the geometric version of Riemann-Roch implies 
$h^0(D(x) + x) = 2$. So the linear system $|D(x) + x|$ is a $g^1_5$. But $D(x)+x$ is part if a fibre of $\pi: X \lra \P^1$
and the corresponding linear system does not admit fixed points. So this cannot occur completing the proof of the proposition. $\hfill{\square}$ 
\end{pf}

\noindent
Combining Propositions 6.1 and 6.3 we proved the theorem stated in the introduction.\\

\noindent
Let $(A, \Theta)$ denote a principally polarized abelian variety of dimension $g$. The cohomology class 
$\frac{1}{(g-1)!}[\Theta]^{g-1}$ is not divisible in $H^{2g-2}(A,\Z)$ therefore called the dimension-one 
{\it minimal cohomology class} of $(A,\Theta)$. If $g \geq 3$ it is not at all clear whether a multiple of it 
contains a smooth irreducible curve. We obtain as a consequence of Theorem 5.3 and the result of this section

\begin{cor}
Let $X_1$ and $X_2$ be trigonal curves of positive genus with simple ramification and disjoint branching and let $\Theta$ denote the 
canonical product principal polarization of $JX_1 \times JX_2$. If $g := g(X_1) + g(X_2) \geq 3$, then the class 
$\frac{3}{(g-1)!}[\Theta]^{g-1}$ is represented by a smooth irreducible curve.
\end{cor} 

\noindent
For the proof we need the following well known lemma, for which we include a proof in lack of a reference:

\begin{lem}
For $i=1,2$ let $(A_i,\Theta_i)$ be principally polarized abelian varieties with $\End_\Q(A_i) = \Q$ and $Hom(A_1,A_2) = 0$. 
Then 
$A_1 \times A_2$ 
admits no principal polarization apart from the canonical product polarization.  
\end{lem}

\begin{pf}
Let $\Theta$ denote the canonical product polarization of $A_1 \times A_2$. For any polarization $L$ let 
$\varphi_L: A_1 \times A_2 \lra (A_1 \times A_2)^*$ denote the associated isogeny onto the dual abelian variety.
According to \cite{lb}, Theorem 5.2.8 the map 
$NS(A_1 \times A_2) \lra \End (A_1 \times A_2), \,\, L \mapsto \varphi_L \circ \varphi_{\Theta}^{-1}$
induces a bijection between the sets of principal polarizations and totally positive automorphisms of 
$A_1 \times A_2$ symmetric with respect to the Rosati involution. But $\End(A_1 \times A_2) = \Z \times \Z$ and the 
only such automorphism is $1_{A_1} \times 1_{A_2}$, the Rosati involution being the identity. Since  
$1_{A_1} \times 1_{A_2}$ corresponds to the polarization $\Theta$, this implies the assertion.$\hfill{\square}$
\end{pf}

\noindent
{\bf Proof} (of Corollary 6.4): By Theorem 5.3 it suffices to show that the principal polarization $\Xi$ of $P$
coincides with the canonical product principal polarization of $JX_1 \times JX_2$. According to Lemma 6.5 it is enough to
show that for general trigonal curves as in the corollary we have $\End_{\Q}(JX_i) = \Q$ for $i=1,2$ and $Hom(JX_1,JX_2) = 0$,
since if the polarizations coincide for general trigonal curves they do so for all such curves.\\
The second condition being obvious it suffices to show that a general trigonal curve $X$ satisfies $\End_{\Q}(JX) = \Q$. 
But for $1 \leq g(X) \leq 4$ any curve is trigonal and for $g(X) \geq 5$ the subspace of trigonal curves in the moduli
space $\M_{g(X)}$ is of dimension $2g(X) + 1$. Hence the assertion follows from the main result of \cite{cgt}
which says that the Jacobian of a general member of a family of curves of genus g and dimension $> 2g-2$ has endomorphism 
algebra $\Q$.  $\hfill{\square}$

\end{document}